\documentclass[a4paper]{amsart}
\usepackage{amsthm,amsfonts,amsmath,amssymb}
\usepackage[abs]{overpic}
\usepackage{comment} 

\newtheorem{theorem}{Theorem}[section]
\newtheorem{lemma}[theorem]{Lemma}

\newtheorem{proposition}[theorem]{Proposition}

\theoremstyle{definition}

\newtheorem{remark}[theorem]{Remark}

\theoremstyle{remark}

\newcommand{\Z}{\mathbb{Z}}

\makeatletter

\@addtoreset{figure}{section}
\makeatother

\makeatletter

\@addtoreset{table}{section}
\makeatother

\makeatletter
  
  \@addtoreset{equation}{section}
\makeatother

\begin{document}
\title{Hurwitz equivalence in the universal dihedral quandle}

\author{Takuji NAKAMURA}
\address{Faculty of Education, 
University of Yamanashi,
Takeda 4-4-37, Kofu, Yamanashi 400-8510, Japan}
\email{takunakamura@yamanashi.ac.jp}

\author{Yasutaka NAKANISHI}
\address{Department of Mathematics, Kobe University, 
Rokkodai-cho 1-1, Nada-ku, Kobe 657-8501, Japan}
\email{nakanisi@math.kobe-u.ac.jp}

\author{Shin SATOH}
\address{Department of Mathematics, Kobe University, 
Rokkodai-cho 1-1, Nada-ku, Kobe 657-8501, Japan}
\email{shin@math.kobe-u.ac.jp}

\author[Kodai Wada]{Kodai Wada}
\address{Department of Mathematics, Kobe University, Rokkodai-cho 1-1, Nada-ku, Kobe 657-8501, Japan}
\email{wada@math.kobe-u.ac.jp}

\makeatletter
\@namedef{subjclassname@2020}{%
  \textup{2020} Mathematics Subject Classification}
\makeatother
\subjclass[2020]{20F36, 05E18, 57K12, 57K10}

\keywords{Hurwitz action, braid, quandle, pure braid, virtual braid}

\thanks{This work was supported by JSPS KAKENHI Grant Numbers 
JP20K03621, JP22K03287, and JP23K12973.}



\begin{abstract}
We investigate the Hurwitz action of the $m$-braid group on the $m$-fold Cartesian product of the universal dihedral quandle.
We introduce three computable invariants and 
prove that they give a complete classification of the orbits under this action. 
As a consequence, we describe an explicit complete system of orbit representatives. 
We further obtain analogous classifications 
for the corresponding Hurwitz actions of the pure $m$-braid group, 
the virtual $m$-braid group, and the virtual pure $m$-braid group. 
\end{abstract}

\maketitle

\section{Introduction} 
For a group $G$, 
the Hurwitz action of the $m$-braid group $B_m$ 
on the $m$-fold Cartesian product $G^m$ 
is a fundamental object in combinatorial group theory 
and topology. 
This action is defined by the standard generator $\sigma_i$ 
($i=1,\dots,m-1$) of $B_m$, which acts as 
\[
(g_1,\dots,g_i,g_{i+1},\dots,g_m)\cdot\sigma_i=
(g_1,\dots,g_{i-1},g_{i+1},g_{i+1}^{-1}g_ig_{i+1},g_{i+2},\dots,g_m)\in G^{m}.
\] 
The problem of determining when two elements of $G^m$ 
lie in the same orbit under the Hurwitz action 
has been extensively studied for various classes of groups, 
including symmetric groups~\cite{BT}, 
generalized quaternion groups~\cite{Hou}, 
dicyclic groups~\cite{Sia}, 
and dihedral groups~\cite{Ber, Hou, Sia}. 

The importance of the Hurwitz action stems from 
its deep connections with low-dimensional topology. 
In particular, the classification of surface braids corresponds to 
the classification of Hurwitz equivalence classes of 
their monodromy systems \cite{Kam01,Kam02}. 
Thus, understanding these orbits is 
essential for characterizing topological invariants of 
branched coverings and surface-knots. 

The Hurwitz action on groups admits a natural extension to quandles. 
A quandle~\cite{Joy,Mat} is an algebraic structure 
equipped with a binary operation $*$, 
whose axioms algebraically encode the Reidemeister moves in knot theory. 
For a quandle $X$, the $m$-braid group $B_m$ acts on $X^m$ 
by 
\[
(x_1,\dots,x_i,x_{i+1},\dots,x_m)\cdot\sigma_i=
(x_1,\dots,x_{i-1},x_{i+1},x_i*x_{i+1},x_{i+2},\dots,x_m)\in X^{m}.
\]
Since any conjugacy class of a group $G$ carries a natural quandle structure 
under the operation $x*y=y^{-1}xy$, 
results on the Hurwitz action for groups 
can be reinterpreted within the framework of quandle theory. 
For example, the set of all reflections 
in the dihedral group of order $2n$ forms 
a conjugacy class isomorphic to 
the order-$n$ dihedral quandle $R_n=({\Z}/n{\Z},*)$, 
where $a*b\equiv 2b-a\pmod{n}$. 
Accordingly, results for the Hurwitz action on 
the dihedral group translate directly into 
statements about the action on $(R_n)^m$. 

The Hurwitz action on $X^m$ also admits an interpretation 
in terms of quandle colorings. 
If $v,w\in X^m$ satisfy that $v\cdot\beta=w$ for some $\beta\in B_m$, 
then there exists an $X$-coloring of a braid diagram of $\beta$ 
whose colors at the top and bottom endpoints are given by $v$ and $w$, respectively. 
Colorings, introduced by Fox \cite{Fox}, 
play a central role in knot invariants. 
In particular, Fox $n$-colorings coincide with  
colorings by the dihedral quandle $R_n$. 
This connection further motivates the study of Hurwitz actions on dihedral quandles. 

In this paper, we focus on the universal dihedral quandle 
$R_{\infty}=({\Z},*)$, 
where $a*b=2b-a$. 
This quandle may be regarded as the universal covering 
of the finite dihedral quandles $R_n$. 
We introduce three invariants associated with 
an element $v\in(R_{\infty})^m$, 
namely, two integers $\Delta(v)$ and $d(v)$, 
together with a multiset $M(v)$, 
and show that they completely determine the Hurwitz orbit. 

\begin{theorem}\label{thm11}
For $v,w\in(R_{\infty})^m$ with $m\geq 2$, 
the following are equivalent. 
\begin{enumerate}
\item
There exists an $m$-braid $\beta\in B_m$ 
such that $v\cdot\beta=w$. 
\item
$\Delta(v)=\Delta(w)$, $d(v)=d(w)$, and $M(v)=M(w)$. 
\end{enumerate}
\end{theorem}

Moreover, we provide a complete system of orbit representatives 
under the Hurwitz action. 
Specifically, we show that every orbit contains a representative 
of the form \[(x,\dots,x,z,y\dots,y),\] 
where $x,y,z\in R_{\infty}$ satisfy certain explicit conditions (Propositions~\ref{prop33}, \ref{prop49}, and~\ref{prop54}). 

We further extend these results in several directions. 
First, we analyze the Hurwitz action of the pure $m$-braid group $P_m$. 
We then turn to the virtual setting, 
investigating  the actions of the virtual $m$-braid group $VB_m$ 
and the virtual pure $m$-braid group $VP_m$, 
which naturally extend the classical braid group $B_{m}$ 
and its pure subgroup $P_m$, respectively. 
In this broader framework, we introduce a refinement $M^*(v)$ 
of the invariant $M(v)$. 
This refined invariant enables us to obtain precise 
characterizations of the orbits under each of these actions, 
as summarized in the following results. 

\begin{theorem}\label{thm12}
For $v,w\in(R_{\infty})^m$ with $m\geq 2$, 
the following are equivalent. 
\begin{enumerate}
\item
There exists a pure $m$-braid $\beta\in P_m$ 
such that $v\cdot\beta=w$. 
\item
$\Delta(v)=\Delta(w)$, $d(v)=d(w)$, and $M^*(v)=M^*(w)$. 
\end{enumerate}
\end{theorem}

\begin{theorem}\label{thm13}
For $v,w\in(R_{\infty})^m$ with $m\geq 2$, 
the following are equivalent. 
\begin{enumerate}
\item
There exists a virtual $m$-braid $\beta\in VB_m$ 
such that $v\cdot\beta=w$. 
\item
$d(v)=d(w)$ and $M(v)=M(w)$. 
\end{enumerate}
\end{theorem}

\begin{theorem}\label{thm14}
For $v,w\in(R_{\infty})^m$ with $m\geq 2$, 
the following are equivalent. 
\begin{enumerate}
\item
There exists a virtual pure $m$-braid $\beta\in VP_m$ 
such that $v\cdot\beta=w$. 
\item
$d(v)=d(w)$ and $M^*(v)=M^*(w)$. 
\end{enumerate}
\end{theorem}

This paper is organized as follows. 
In Section~\ref{sec2}, we introduce the three quantities 
$\Delta(v)$, $d(v)$, and $M(v)$ for an element $v\in {\Z}^m$, 
and prove their invariance under the Hurwitz action (Lemma~\ref{lem22}). 
Sections~\ref{sec3}--\ref{sec5} are devoted to 
the proof of Theorem~\ref{thm11}, 
covering the cases $m=2$ (Theorem~\ref{thm31}), 
odd $m\geq 3$ (Theorem~\ref{thm47}), and even $m\geq 4$ 
(Theorem~\ref{thm53}), respectively. 
Within these sections, 
we also establish a complete system of orbit representatives. 
Section~\ref{sec6} examines the action by the pure $m$-braid group $P_m$, 
introducing the refined invariant $M^*(v)$ 
and providing a characterization of its orbit (Theorem~\ref{thm63}). 
Finally, in Section~\ref{sec7}, 
we investigate the actions of the virtual $m$-braid group $VB_m$ and 
the virtual pure $m$-braid group $VP_m$, 
and establish corresponding orbit characterizations 
(Theorem~\ref{thm76}). 

For the sake of simplicity, throughout the rest of this paper, 
we identify the universal dihedral quandle $R_{\infty}$ 
with the set of integers ${\Z}$.

\section{Preliminaries}\label{sec2}

\begin{figure}[htbp]
  \vspace{1em}
  \begin{center}
   \hspace{1em}
    \begin{overpic}[]{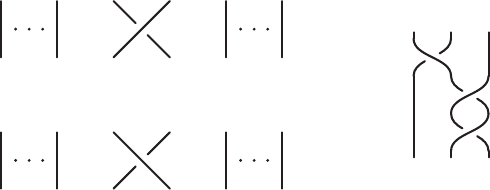}
      \put(-26,75){$\sigma_{i}=$}
      \put(-34,11.5){$\sigma_{i}^{-1}=$}
      \put(-4,98){$a_{1}$}
      \put(18,98){$a_{i-1}$}
      \put(47,98){$a_{i}$}
      \put(74,98){$a_{i+1}$}
      \put(102,98){$a_{i+2}$}
      \put(131,98){$a_{m}$}
      \put(-4,54){$a_{1}$}
      \put(15,54){$a_{i-1}$}
      \put(37,54){$a_{i+1}$}
      \put(58,54){$2a_{i+1}-a_{i}$}
      \put(106,54){$a_{i+2}$}
      \put(131,54){$a_{m}$}
      \put(-4,34){$a_{1}$}
      \put(18,34){$a_{i-1}$}
      \put(47,34){$a_{i}$}
      \put(74,34){$a_{i+1}$}
      \put(102,34){$a_{i+2}$}
      \put(131,34){$a_{m}$}
      \put(-4,-11){$a_{1}$}
      \put(15,-11){$a_{i-1}$}
      \put(39,-11){$2a_{i}-a_{i+1}$}
      \put(88,-11){$a_{i}$}
      \put(102,-11){$a_{i+2}$}
      \put(131,-11){$a_{m}$}
      \put(243,43){$=\sigma_{1}^{-1}\sigma_{2}^{2}$}
      \put(197,80){$1$}
      \put(210,80){$-5$}
      \put(232,80){$4$}
      \put(196.5,5){$7$}
      \put(215,5){$7$}
      \put(230,5){$10$}
    \end{overpic}
  \end{center}
  \vspace{1em}
  \caption{The Hurwitz action of $B_{m}$ on $\Z^{m}$}
  \label{Hurwitz-action}
\end{figure}

For an integer $m\geq2$, 
let $B_m$ denote the $m$-braid group 
with the standard generators 
$\sigma_1,\dots,\sigma_{m-1}$. 
The set of $m$-tuples of integers, 
\[\Z^m=\{(a_1,\dots,a_m)\mid a_1,\dots,a_m\in\Z\},\]
admits the {\it Hurwitz action} of $B_m$ defined by 
\begin{align*}
&(a_{1},\dots,a_{i},a_{i+1},\dots,a_{m})\cdot\sigma_i  =(a_1,\dots,a_{i-1}, a_{i+1},2a_{i+1}-a_i,a_{i+2},\dots,a_m), \\
&(a_{1},\dots,a_{i},a_{i+1},\dots,a_{m})\cdot\sigma_i^{-1}=(a_1,\dots,a_{i-1}, 
2a_i-a_{i+1},a_i,a_{i+2},\dots,a_m). 
\end{align*} 
See the left of Figure~\ref{Hurwitz-action}. 
We say that two elements $v,w\in\Z^m$ are {\it equivalent} 
if there exists an $m$-braid $\beta\in B_m$ 
such that $v\cdot\beta=w$, 
and denote it by $v\sim w$. 
For example, the right of the figure shows that 
\[(1,-5,4)\cdot \sigma_1^{-1}\sigma_2^2=(7,7,10)\in\Z^{3},\]
which means that $(1,-5,4)\sim(7,7,10)$.

An element $v\in\Z^m$ is called {\it trivial} 
if $v=(a,\dots,a)=a\cdot {\mathbf 1}$ for some $a\in\Z$, 
where ${\mathbf 1}=(1,\dots,1)$. 
By definition, if $v\sim w$ and $v$ is trivial, 
then $v=w$; 
that is, the orbit of a trivial element $v$ 
consists solely of $v$. 

Now, we introduce three quantities for 
an element $v=(a_1,\dots,a_m)\in{\Z}^m$. 
The first quantity $\Delta(v)$ is an integer defined by 
\[\Delta(v)=\sum_{i=1}^m(-1)^{i-1} a_i.\] 
The second quantity $d(v)$ is a nonnegative integer defined by 
\[
d(v)=\begin{cases}
0 & \text{if $v$ is trivial}, \\
\gcd\{a_i-a_j\mid 1\leq i,j\leq m\}>0 
& \text{if $v$ is nontrivial}. 
\end{cases}
\]
We remark that 
$d(v)=\gcd\{a_i-a_1\mid 1\leq i\leq m\}$ 
for $v$ nontrivial. 
The third quantity $M(v)$ is defined as  
an unordered $m$-tuple 
consisting of the congruence classes of 
$a_1,\dots,a_m$ modulo $2d(v)$;  
that is, $M(v)$ is the multiset 
\[M(v)=\{a_1,\dots,a_m\ ({\rm mod}~2d(v))\}.\]
Here, we use the convention that if $d(v)=0$, 
then $a_i$ (mod~$0$) simply denotes the integer $a_i\in{\Z}$ itself. 
For example, the element $v=(4,1,7,-2,-5)\in{\Z}^5$ 
satisfies that   
\[\Delta(v)=7, \ d(v)=3, \mbox{ and } 
M(v)=\{1,1,1,4,4\}.\]

\begin{lemma}\label{lem21}
For $v=(a_1,\dots,a_m)\in\Z^m$, 
the following hold. 
\begin{enumerate}
\item
$a_1\equiv\dots\equiv a_m \pmod{d(v)}$. 
\item
If $m$ is odd, then 
$\Delta(v)\equiv a_1\pmod{d(v)}$. 
\item
If $m$ is even, 
then $\Delta(v)\equiv 0\pmod{d(v)}$. 
\item
If $v$ is nontrivial, 
then 
\[M(v)=\{\underbrace{r,\dots,r}_{p},
\underbrace{r+d(v),\dots,r+d(v)}_{m-p}\}\]
for some $0\leq r<d(v)$ and 
$1\leq p\leq m-1$. 
\end{enumerate}
\end{lemma}

\begin{proof}
(i) Since $a_i-a_1$ is divisible by $d(v)$, 
we have $a_i\equiv a_1\pmod{d(v)}$ for any $i$. 

(ii) It follows from (i) that 
$\Delta(v)=a_1-\sum_{i=1}^{(m-1)/2}(a_{2i}-a_{2i+1})
\equiv a_1\pmod{d(v)}$. 

(iii)
It follows from (i) that 
$\Delta(v)=\sum_{i=1}^{m/2}(a_{2i-1}-a_{2i})\equiv 0\pmod{d(v)}$. 

(iv) This is a direct consequence of (i). 
\end{proof}

\begin{lemma}\label{lem22}
If $v\sim w$, then the following hold. 
\begin{enumerate}
\item
$\Delta(v)=\Delta(w)$.  
\item
$d(v)=d(w)$. 
\item
$M(v)=M(w)$.
\end{enumerate}
\end{lemma}

\begin{proof}
Let $v=(a_1,\dots,a_m)$. 
It suffices to consider the case 
$w=v\cdot \sigma_i$ ($i=1,\dots,m-1$). 

(i) We have
\[\begin{split}
\Delta(w)-\Delta(v)
=&[(-1)^{i-1}a_{i+1}+(-1)^i(2a_{i+1}-a_i)]\\
&-[(-1)^{i-1}a_i+(-1)^ia_{i+1}]=0.
\end{split}\]

(ii) For $i=1$, we have 
\[\begin{split}
d(w)&=
\gcd\{(2a_2-a_1)-a_2,a_3-a_2,\dots,a_m-a_2\}\\
&=\gcd\{a_2-a_1,(a_3-a_1)-(a_2-a_1),\dots,(a_m-a_1)-(a_2-a_1)\}\\
&=\gcd\{a_2-a_1,a_3-a_1,\dots,a_m-a_1\}=d(v). 
\end{split}\]

For $i=2,\dots,m-1$, we have
\[\begin{split}
d(w)&=
\gcd\{a_2-a_1,\dots,a_{i-1}-a_1,\\
& \hspace{50pt}
a_{i+1}-a_1,(2a_{i+1}-a_i)-a_1,
a_{i+2}-a_1,\dots,a_m-a_1\}\\
&=
\gcd\{a_2-a_1,\dots,a_{i-1}-a_1,\\
& \hspace{50pt}
a_{i+1}-a_1,2(a_{i+1}-a_1)-(a_i-a_1),
a_{i+2}-a_1,\dots,a_m-a_1\}\\
&=
\gcd\{a_2-a_1,\dots,a_m-a_1\}=d(v). 
\end{split}\]

(iii) By Lemma~\ref{lem21}(i), we have
\[2a_{i+1}-a_i=2(a_{i+1}-a_i)+a_i\equiv a_i ~({\rm mod}~2d(v)).\] 
Therefore, it follows from (ii) that 
\begin{align*}
M(w)&=\{\dots,a_{i+1},2a_{i+1}-a_i,\dots \ ({\rm mod}~2d(w))\}\\
&=\{\dots,a_i,a_{i+1},\dots \ ({\rm mod}~2d(v))\}=M(v). 
\end{align*}
\end{proof}

Let $S_m$ denote the symmetric group on $\{1,\dots,m\}$. 
For an $m$-braid $\beta\in B_m$, 
we denote by $\pi_{\beta}\in S_m$ 
the permutation associated with $\beta$; 
that is, $\beta$ connects each $k$th top point from the left 
with the $\pi_{\beta}(k)$th bottom point ($k=1,\dots,m$).

\begin{lemma}\label{lem23}
Suppose that $v=(a_1,\dots,a_m)$,  $w=(b_1,\dots,b_m)\in\Z^m$, 
and $\beta\in B_m$ satisfy $v\cdot \beta=w$. 
Then 
\[b_{\pi_\beta(k)}\equiv a_k \pmod{2d(v)}\]  
for any $k=1,\dots,m$. 
\end{lemma}

\begin{proof}
It suffices to consider the case 
$\beta=\sigma_i$ $(i=1,\dots,m-1)$. 
For $k\ne i,i+1$, 
we have $b_{\pi_{\beta}(k)}=b_k=a_k$. 
For $k=i$, we have 
\[b_{\pi_{\beta}(i)}=b_{i+1}
=2a_{i+1}-a_i
\equiv a_i\pmod{2d(v)}.\]
For $k=i+1$, we have 
$b_{\pi_{\beta}(i+1)}=b_i
=a_{i+1}$. 
\end{proof}

\section{The case $m=2$}\label{sec3}

In this section, 
we study the equivalence relation $\sim$ on $\Z^2$. 

\begin{theorem}\label{thm31}
For $v,w\in\Z^2$, 
the following are equivalent. 
\begin{enumerate}
\item
$v\sim w$. 
\item
$\Delta(v)=\Delta(w)$, $d(v)=d(w)$, and $M(v)=M(w)$. 
\end{enumerate}
\end{theorem}

\begin{proof}
\underline{(i)$\Rightarrow$(ii).} 
This follows from Lemma~\ref{lem22}. 

\underline{(ii)$\Rightarrow$(i).} 
Let $\Delta=\Delta(v)=\Delta(w)$. 
Then we have 
\[v\sim(r,r-\Delta) \mbox{ and }w\sim(s,s-\Delta)\] 
for some $0\leq r,s\leq |\Delta|$. 
In fact, we have 
\[v\sim v\cdot \sigma_1^n=v-n\Delta \cdot{\mathbf 1}\]
for any $n\in{\Z}$. 
Since $d(v)=d(w)=|\Delta|$, we have 
\[M(v)=\{r,r+|\Delta|\}\mbox{ and }M(w)=\{s,s+|\Delta|\}.\] 
Moreover, since $M(v)=M(w)$, we have $r=s$, and hence, 
$v\sim w$. 
\end{proof}

\begin{remark}\label{rem32}
For $v\in{\Z}^2$, we have $d(v)=|\Delta(v)|$ 
as seen in the proof of Theorem~\ref{thm31}. 
Therefore, $v\sim w\in{\Z}^2$ if and only if 
$\Delta(v)=\Delta(w)$ and $M(v)=M(w)$. 
\end{remark}

\begin{proposition}\label{prop33}
A complete system of representatives 
of $\Z^2/\sim$ is given by 
\[\{(x,x)\mid x\in{\Z}\}\sqcup
\{(x,y)\mid 0\leq 2x<y\}\sqcup
\{(x,y)\mid 0\leq 2y<x\}.\] 
\end{proposition}

\begin{proof}
By the proof of Theorem~\ref{thm31}, 
any nontrivial element in ${\Z}^2$ is equivalent to 
$(r,r-\Delta)\sim(r+\Delta,r)$ for some $\Delta\ne 0$ and $0\leq r< |\Delta|$. 
\begin{itemize}
\item
If $\Delta<0$, then we have $0\leq 2r<r-\Delta$, 
and hence, $(r,r-\Delta)$ belongs to the second set 
in the system. 
\item
If $\Delta>0$, then we have $0\leq 2r<r+\Delta$, 
and hence, $(r+\Delta,r)$ belongs to the third set 
in the system above. 
\end{itemize}

For two nontrivial elements 
$v_i=(x_i,y_i)$ ($i=1,2$) in the system, 
we have  
\[\Delta(v_i)=x_i-y_i, \ d(v_i)=|y_i-x_i|, \mbox{ and }
M(v_i)=\{x_i,y_i\}.\] 
If $\Delta(v_i)<0$, then $v_i$ belongs to the second set. 
Since we have 
\[d(v_i)=y_i-x_i>0 \mbox{ and }0\leq 2x_i<y_i,\]
it follows that 
\[0\leq x_i<d(v_i) \mbox{ and }d(v_i)\leq y_i<2d(v_i).\]
Therefore, if $v_1\sim v_2$, 
then we have $x_1=x_2$ and $y_1=y_2$ by Lemma~\ref{lem22}, 
and hence, $v_1=v_2$. 
The case $\Delta(v_i)>0$ is treated analogously. 
\end{proof}

\section{The case $m$ odd}\label{sec4}

In this section, 
we study the equivalence relation $\sim$ on $\Z^m$ 
for $m=2k-1\geq 3$. 

For an element $v=(a_1,a_2,a_3)\in\Z^3$, 
we define a nonnegative integer by  
\[|v|=\max\{a_1,a_2,a_3\}-\min\{a_1,a_2,a_3\}\geq0.\] 
We remark $|v|=0$ if and only if $v$ is trivial. 

\begin{lemma}\label{lem41}
For $v=(a_1,a_2,a_3)\in{\Z}^3$, 
if $a_1$, $a_2$, and $a_3$ are mutually distinct, 
then $v$ is equivalent to some $w\in\Z^3$ 
with $|w|<|v|$. 
\end{lemma}

\begin{figure}[htbp]
  \centering
    \vspace{3em}
    \begin{overpic}[width=9cm]{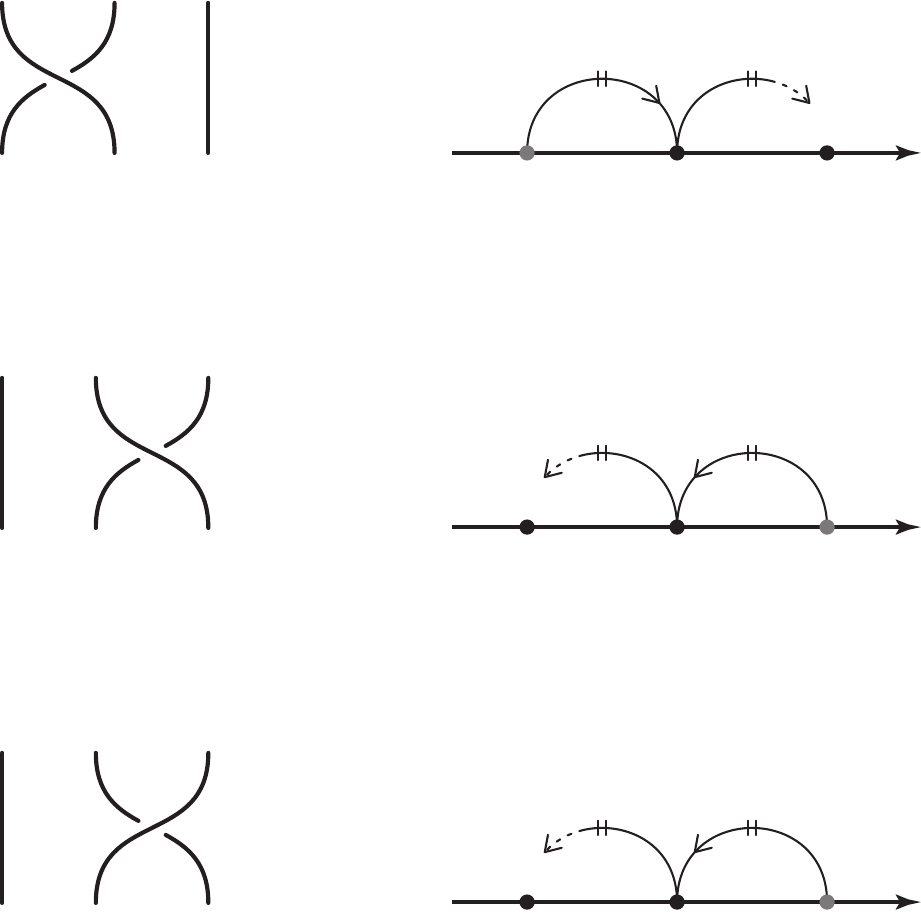}
       \put(-15,280){\underline{$a_{2}<a_{1}$}}
      \put(-4,262){$a_{1}$}
      \put(27,262){$a_{2}$}
      \put(53,262){$a_{3}$}
      \put(-20.5,202){$2a_{1}-a_{2}$}
      \put(27,202){$a_{1}$}
      \put(53,202){$a_{3}$}
      \put(142,203){$a_{2}$}
      \put(184,203){$a_{1}$}
      \put(226,203){$a_{3}$}
      \put(219,251.5){$a_{1}+(a_{1}-a_{2})$}
      \put(248,239){\rotatebox[origin=c]{90}{$=$}}
      \put(228,224){\framebox{$2a_{1}-a_{2}$}}
       \put(-15,175){\underline{$a_{1}<a_{2}<a_{3}$}}
      \put(-4,157){$a_{1}$}
      \put(22,157){$a_{2}$}
      \put(53,157){$a_{3}$}
      \put(-4,98){$a_{1}$}
      \put(12,98){$2a_{2}-a_{3}$}
      \put(53,98){$a_{2}$}
      \put(142,98){$a_{1}$}
      \put(184,98){$a_{2}$}
      \put(226,98){$a_{3}$}
      \put(96,148.5){$a_{2}-(a_{3}-a_{2})$}
      \put(125,136){\rotatebox[origin=c]{90}{$=$}}
      \put(105,121){\framebox{$2a_{2}-a_{3}$}}
      \put(-15,70.5){\underline{$a_{2}>a_{3}$}}
      \put(-4,52.5){$a_{1}$}
      \put(22,52.5){$a_{2}$}
      \put(53,52.5){$a_{3}$}
      \put(-4,-6){$a_{1}$}
      \put(22,-6){$a_{3}$}
      \put(37.3,-6){$2a_{3}-a_{2}$}
      \put(142,-6){$a_{1}$}
      \put(184,-6){$a_{3}$}
      \put(226,-6){$a_{2}$}
      \put(96,44){$a_{3}-(a_{2}-a_{3})$}
      \put(125,31.5){\rotatebox[origin=c]{90}{$=$}}
      \put(105,16.5){\framebox{$2a_{3}-a_{2}$}}
    \end{overpic}
  \vspace{1em}
  \caption{Proof of Lemma~\ref{lem41}}
  \label{numberline}
\end{figure}

\begin{proof}
We consider only the case $a_1<a_3$; 
the case $a_{1}>a_{3}$ is treated analogously. 
We define $w\in\Z^{3}$ by 
\[w=
\begin{cases}
v\cdot\sigma_1^{-1} & \text{for }a_2<a_1, \\
v\cdot\sigma_2^{-1} & \text{for }a_1<a_2<a_3, \\
v\cdot\sigma_2 & \text{for }a_2>a_3.
\end{cases}\]
Then we have $|w|<|v|$. 
In fact, for the first case, we have 
\[\begin{split}
|w|=|(2a_1-a_2,a_1,a_3)|
&=\max\{2a_1-a_2,a_3\}-a_1\\
&=\max\{a_1-a_2,a_3-a_1\}<a_3-a_2=|v|.
\end{split}
\]
See the top of Figure~\ref{numberline}. 
For the second case, we have 
\[
\begin{split}
|w|=|(a_1,2a_2-a_3,a_2)|
&=a_2-\min\{a_1,2a_2-a_3\}\\
&=\max\{a_2-a_1,a_3-a_2\}<a_3-a_1=|v|.
\end{split}
\]
See the middle of the figure. 
For the third case, we have 
\[
\begin{split}
|w|=|(a_1,a_3,2a_3-a_2)|
&=a_3-\min\{a_1,2a_3-a_2\}\\
&=\max\{a_3-a_1,a_2-a_3\}<a_2-a_1=|v|.
\end{split}
\]
See the bottom of the figure. 
\end{proof}

\begin{lemma}\label{lem42}
Any element $v\in\Z^3$ 
is equivalent to 
$(x,y,y)$ for some  $x\leq y$. 
\end{lemma}

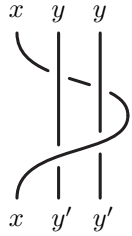
\begin{figure}[htb]
  \centering
    \vspace{1em}
    \begin{overpic}[]{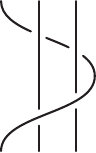}
      \put(-2.5,78){$x$}
      \put(15.5,78){$y$}
      \put(34,78){$y$}
      \put(-2.5,-10){$x$}
      \put(15.5,-10){$y'$}
      \put(34,-10){$y'$}
    \end{overpic}
  \vspace{1em}
  \caption{The braid $\sigma_1\sigma_2^2\sigma_1$}
  \label{sigma1221}
\end{figure}

\begin{proof}
By Lemma~\ref{lem41}, 
there exists a finite sequence of elements in $\Z^3$, 
\[v=v_1,\dots,v_n,\]
such that 
\begin{enumerate}
\item
$v_i=(x_i,y_i,z_i)$ $(i=1,\dots,n)$, 
\item
$v_i\sim v_{i+1}$ and $|v_i|>|v_{i+1}|$ $(i=1,\dots,n-1)$, 
\item
$x_i$, $y_i$, and $z_i$ are mutually distinct $(i=1,\dots,n-1)$, and  
\item
at least two of $x_n$, $y_n$, and $z_n$ are the same. 
\end{enumerate}
If $x_n=z_n$, that is, $v_n=(x_n,y_n,x_n)$, 
then we have 
\[v\sim v_n\cdot\sigma_1^{-1}=(2x_n-y_n,x_n,x_n).\]
If $x_n=y_n$, that is, $v_n=(x_n,x_n,z_n)$, 
then we have 
\[v\sim v_n\cdot\sigma_2^{-1}\sigma_1^{-1}=(z_n,x_n,x_n).\]
Therefore, we may assume that 
$v\sim(x,y,y)$ for some $x,y\in\Z$. 

For $x\leq y$, the conclusion already follows. 
If $x>y$, then we have 
\[
v\sim (x,y,y)\sim 
(x,y,y)\cdot\sigma_1\sigma_2^2\sigma_1
=(x,y',y'),\]
where $y'=2x-y>x$. 
See Figure~\ref{sigma1221}. 
\end{proof}

\begin{lemma}\label{lem43}
Let $v\in\Z^{2k-1}$ be a nontrivial element of the form 
\[v=(\underbrace{x,\dots,x}_{2p-1},
\underbrace{y,\dots,y}_{2q},
\underbrace{z,\dots,z}_{2r})\]
for some $x,y,z,p,q,r\in\Z$ with $p,q,r\geq 1$ and $p+q+r=k$. 
If $x<y<z$, 
then $v$ is equivalent to 
\[w=
(\underbrace{x,\dots,x}_{2p-1},
\underbrace{y',\dots,y'}_{2q'}, 
\underbrace{z',\dots,z'}_{2r'})\]
for some $y', z', q', r'\in\Z$ with 
$x\leq y'\leq z'<z$ and $\{q',r'\}=\{q,r\}$. 
\end{lemma}

\begin{figure}[htbp]
\vspace{1em}
  \begin{flushright}
    \begin{overpic}[width=12cm]{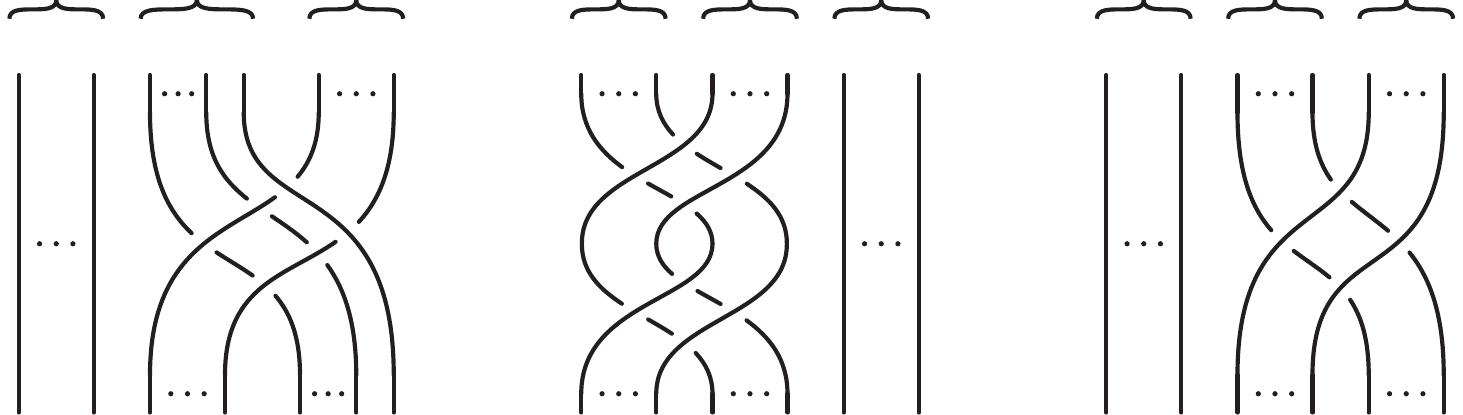}
      \put(-17,37.5){$\beta=$}
      \put(-0.5,101){{\footnotesize $2p-1$}}
      \put(41.5,101){{\footnotesize $2q$}}
      \put(79,101){{\footnotesize $2r$}}
      \put(1.5,83){$x$}
      \put(19,83){$x$}
      \put(32.5,84){$y$}
      \put(45.5,84){$y$}
      \put(55,84){$y$}
      \put(72,83){$z$}
      \put(89.5,83){$z$}
      \put(1.5,-10){$x$}
      \put(19,-10){$x$}
      \put(32,-10){$z'$}
      \put(49.5,-10){$z'$}
      \put(67,-10){$y$}
      \put(80,-10){$y$}
      \put(89.5,-10){$y$}

      \put(112,37.5){$\beta'=$}
      \put(130.5,101){{\footnotesize $2p-1$}}
      \put(170.5,101){{\footnotesize $2r$}}
      \put(201,101){{\footnotesize $2q$}}
      \put(133,83){$x$}
      \put(150.5,83){$x$}
      \put(164,83){$z'$}
      \put(180.5,83){$z'$}
      \put(194.5,84){$y$}
      \put(212,84){$y$}
      \put(133,-10){$x$}
      \put(150.5,-10){$x$}
      \put(163.5,-10){$z''$}
      \put(180,-10){$z''$}
      \put(194.5,-10){$y$}
      \put(212,-10){$y$}

      \put(232,37.5){$\beta''=$}
      \put(253.5,100){{\footnotesize $2p-1$}}
      \put(293.3,101){{\footnotesize $2r$}}
      \put(323.5,102){{\footnotesize $2q$}}
      \put(255,83){$x$}
      \put(273,83){$x$}
      \put(286.5,83){$z''$}
      \put(302.5,83){$z''$}
      \put(317,84){$y$}
      \put(335,84){$y$}
      \put(255,-10){$x$}
      \put(273,-10){$x$}
      \put(286.5,-10){$y$}
      \put(303.5,-10){$y$}
      \put(316,-10){$z''$}
      \put(334,-10){$z''$}
    \end{overpic}
  \end{flushright}
  \vspace{1em}
  \caption{Three braids $\beta$, $\beta'$, and $\beta''$}
  \label{betas}
\end{figure}

\begin{proof}
Let $\beta,\beta',\beta''\in B_{2k-1}$ be braids  
as shown in Figure~\ref{betas}, 
where $z'=2y-z$ and $z''=2x-z'$. 
Then we define $w\in\Z^{2k-1}$ by 
\[w=
\begin{cases}
v\cdot\beta & \text{for }y<z\leq y+(y-x), \\
v\cdot\beta\beta' & \text{for }y+(y-x)<z\leq y+2(y-x), \\
v\cdot\beta\beta'\beta'' & \text{for }z>y+2(y-x).
\end{cases}\]
For the first case, we have 
$x\leq z'<y<z$ and 
\[
v\sim w=
(\underbrace{x,\dots,x}_{2p-1},
\underbrace{z',\dots,z'}_{2r},
\underbrace{y,\dots,y}_{2q}). 
\]
For the second case, we have 
$x< z''\leq y<z$ and 
\[v\sim w=
(\underbrace{x,\dots,x}_{2p-1},
\underbrace{z'',\dots,z''}_{2r},
\underbrace{y,\dots,y}_{2q}).
\]
For the third case, 
we have $x<y<z''<z$ and 
\[v\sim w=
(\underbrace{x,\dots,x}_{2p-1},
\underbrace{y,\dots,y}_{2q},
\underbrace{z'',\dots,z''}_{2r}).
\]
\end{proof}

\begin{lemma}\label{lem44}
Let $v\in\Z^{2k-1}$ be a nontrivial element of the form 
\[
v=(\underbrace{x,\dots,x}_{2p-1},
\underbrace{y,\dots,y}_{2q},
\underbrace{z,\dots,z}_{2r})\]
for some $x,y,z,p,q,r\in\Z$ 
with $p,q,r\geq 1$ and $p+q+r=k$. 
If $x<y<z$, then 
$v$ is equivalent to 
\[
w=(\underbrace{x,\dots,x}_{2p'-1},
\underbrace{y',\dots,y'}_{2k-2p'})\]
for some $y', p'\in\Z$ 
with $x<y'$ and $p'\in\{p,p+q,p+r\}$. 
\end{lemma}

\begin{proof}
By Lemma~\ref{lem43}, 
there exists a finite sequence of elements in $\Z^{2k-1}$, 
\[v=v_1,v_2,\dots,v_n=w,\]
such that 
\begin{enumerate}
\item
$\displaystyle{
v_i=(\underbrace{x,\dots,x}_{2p-1}, 
\underbrace{y_i,\dots,y_i}_{2q_i}, 
\underbrace{z_i,\dots,z_i}_{2r_i})}$ 
$(i=1,\dots,n)$, 
\item
$v_i\sim v_{i+1}$ $(i=1,\dots,n-1)$, 
\item 
$\{q_i,r_i\}=\{q_{i+1},r_{i+1}\}$ $(i=1,\dots,n-1)$, 
\item
$x< y_{i+1}<z_{i+1}<z_i$ $(i=1,\dots,n-2)$, and 
\item 
$x=y_n<z_n$ or $x<y_n=z_n$.
\end{enumerate} 
The conclusion follows from (iii) and (v). 
\end{proof}

\begin{lemma}\label{lem45}
Any element $v\in\Z^{2k-1}$ is 
equivalent to 
\[
(\underbrace{x,\dots,x}_{2p-1},
\underbrace{y,\dots,y}_{2k-2p})\]
for some $x<y$ and $1\leq p\leq k$. 
\end{lemma}

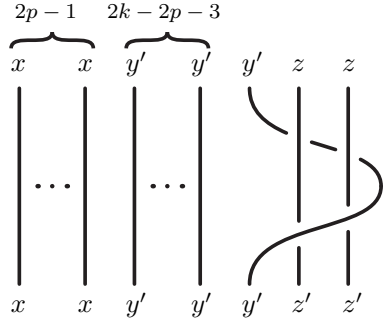
\begin{figure}[htbp]
\vspace{1em}
  \centering
    \begin{overpic}[]{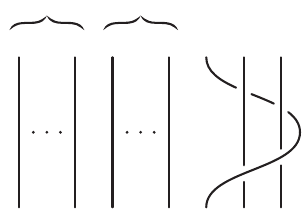}
      \put(10,97){{\footnotesize $2p-1$}}
      \put(45,97){{\footnotesize $2k-2p-3$}}
      \put(6,78){$x$}
      \put(34,78){$x$}
      \put(51,78){$y'$}
      \put(77,78){$y'$}
      \put(97,78){$y'$}
      \put(115,78){$z$}
      \put(133,78){$z$}
      \put(6,-10){$x$}
      \put(34,-10){$x$}
      \put(51,-10){$y'$}
      \put(77,-10){$y'$}
      \put(97,-10){$y'$}
      \put(115,-10){$z'$}
      \put(133,-10){$z'$}
    \end{overpic}
  \vspace{1em}
  \caption{The braid $\sigma_{2k-3}\sigma_{2k-2}^2\sigma_{2k-3}$}
  \label{sigma2k-3}
\end{figure}

\begin{proof}
We proceed by induction on $k\geq 2$. 
The base case $k=2$ follows from Lemma~\ref{lem42}. 

Assume that $k\geq 3$. 
Let $v=(a_1,a_2,\dots,a_{2k-1})$. 
By applying Lemma~\ref{lem42} to 
the last three entries 
$a_{2k-3}$, $a_{2k-2}$, and  $a_{2k-1}$ in $v$, 
we have 
\[v\sim
(a_1,\dots,a_{2k-4},y,z,z)\] 
for some $y\leq z$. 
Moreover, 
by applying the induction hypothesis 
to the first $2k-3$ entries, 
we have 
\[v\sim
(\underbrace{x,\dots,x}_{2p-1},
\underbrace{y',\dots,y'}_{2k-2p-2}, 
z,z)\] 
for some $x< y'$ and $1\leq p\leq k-1$. 

The case $p=k-1$ is straightforward. 
We therefore assume that $1\leq p\leq k-2$. 
If $y'>z$, then we have 
\[v\sim
(\underbrace{x,\dots,x}_{2p-1},
\underbrace{y',\dots,y'}_{2k-2p-2},z,z)
\cdot \sigma_{2k-3}\sigma_{2k-2}^2\sigma_{2k-3} 
=(\underbrace{x,\dots,x}_{2p-1},
\underbrace{y',\dots,y'}_{2k-2p-2},z',z'),\] 
where $z'=2y'-z>y'$. 
See Figure~\ref{sigma2k-3}. 
Hence, we may assume that 
$x< y'\leq z$. 
If $y'=z$, the conclusion already follows. 
If $y'<z$, the conclusion follows from Lemma~\ref{lem44}.
\end{proof}

\begin{proposition}\label{prop46}
Let $v$ be an element in $\Z^{2k-1}$, 
$\Delta=\Delta(v)$, and $d=d(v)$. 
Then there exists an integer $p$ with $1\leq p\leq k$ such that 
\begin{enumerate}
\item
$\displaystyle{
v\sim
(\underbrace{\Delta,\dots,\Delta}_{2p-1}, 
\underbrace{\Delta+d,\dots,\Delta+d}_{2k-2p})}$ and 
\item
$\displaystyle{
M(v)=\{\underbrace{\Delta,\dots,\Delta}_{2p-1},
\underbrace{\Delta+d,\dots,\Delta+d}_{2k-2p}\}}$.
\end{enumerate}
\end{proposition}

\begin{proof}
(i) By Lemma~\ref{lem45}, 
we have 
\[
v\sim (\underbrace{x,\dots,x}_{2p-1},
\underbrace{y,\dots,y}_{2k-2p})\]
for some $x< y$ and $1\leq p\leq k$. 
It follows from Lemma~\ref{lem22}(i) and (ii) that 
$\Delta=x$ and $d=y-x$,  
and hence, $x=\Delta$ and $y=\Delta+d$. 

(ii) This follows from (i) and Lemma~\ref{lem22}(iii). 
\end{proof}

\begin{theorem}\label{thm47}
For $v,w\in\Z^{2k-1}$, 
the following are equivalent. 
\begin{enumerate}
\item 
$v\sim w$. 
\item 
$\Delta(v)=\Delta(w)$, 
$d(v)=d(w)$, and 
$M(v)=M(w)$. 
\end{enumerate}
\end{theorem}

\begin{proof}
\underline{(i)$\Rightarrow$(ii).} 
This follows from Lemma~\ref{lem22}. 

\underline{(ii)$\Rightarrow$(i).} 
Let $\Delta=\Delta(v)=\Delta(w)$ 
and $d=d(v)=d(w)$. 
By Proposition~\ref{prop46}(i), there exist integers $p$ and $q$ 
with $1\leq p,q\leq k$ such that 
\[v\sim
(\underbrace{\Delta,\dots,\Delta}_{2p-1}, 
\underbrace{\Delta+d,\dots,\Delta+d}_{2k-2p})
\mbox{ and }
w\sim
(\underbrace{\Delta,\dots,\Delta}_{2q-1}, 
\underbrace{\Delta+d,\dots,\Delta+d}_{2k-2q}).\]
Since $M(v)=M(w)$ , 
we have $p=q$, and hence, $v\sim w$. 
\end{proof}

\begin{remark}\label{rem48}
For $m=3$, 
the multiset $M(v)$ is uniquely determined by 
$\Delta(v)$ and $d(v)$. 
In fact, we have 
\[M(v)=\{\Delta(v),\Delta(v)+d(v),\Delta(v)+d(v)\},\] 
which includes the case that $v$ is trivial with $d(v)=0$. 
Therefore, $v\sim w\in\Z^{3}$ if and only if 
$\Delta(v)=\Delta(w)$ and $d(v)=d(w)$. 
\end{remark}

\begin{proposition}\label{prop49}
A complete system of representatives 
of $\Z^{2k-1}/\sim$ is given by 
\[\{(x,\dots,x)\in{\Z}^{2k-1}\mid 
x\in{\Z}\}\sqcup \bigsqcup_{p=1}^{k-1}C_p,\] 
where $C_p$ $(1\leq p\leq k-1)$ is the subset of ${\Z}^{2k-1}$ defined by 
\[C_p=\left\{\left.
(\underbrace{x,\dots,x}_{2p-1},
\underbrace{y,\dots,y}_{2k-2p})
\in{\Z}^{2k-1}
\right|\ x<y
\right\}.\] 
\end{proposition}

\begin{proof}
By Proposition~\ref{prop46}(i), 
any nontrivial element in ${\Z}^{2k-1}$ 
is equivalent to 
\[(\underbrace{\Delta,\dots,\Delta}_{2p-1},
\underbrace{\Delta+d,\dots,\Delta+d}_{2k-2p})
\] 
for some $\Delta,d,p\in{\Z}$ with $d> 0$ and $1\leq p\leq k-1$, 
which belongs to the subset $C_{p}$. 

For two nontrivial elements 
\[v_i=(\underbrace{x_i,\dots,x_i}_{2p_i-1},
\underbrace{y_i,\dots,y_i}_{2k-2p_i})\in \bigsqcup_{p
=1}^{k-1}C_{p} \ (i=1,2),\]
we have 
\[\Delta(v_i)=x_i, \ d(v_i)=y_i-x_i, \mbox{ and }
M(v_i)=\{\underbrace{x_i,\dots,x_i}_{2p_i-1},
\underbrace{y_i,\dots,y_i}_{2k-2p_i}\}.\] 
Therefore, if $v_1\sim v_2$, 
then it follows from Lemma~\ref{lem22} that $x_1=x_2$, $y_1=y_2$, and 
$p_1=p_2$, 
and hence, $v_1=v_2$. 
\end{proof}

\section{The case $m\geq 4$ even}\label{sec5}

In this section, 
we study the equivalence relation $\sim$ on $\Z^m$ 
for $m=2k\geq 4$.

\begin{lemma}\label{lem51}
Let $v\in\Z^{2k}$ be a nontrivial element of the form 
\[v=(
\underbrace{x,\dots,x}_{p}, x+\lambda d, 
\underbrace{x+d,\dots,x+d}_{2k-p-1})\] 
for some $x,\lambda,d,p\in\Z$ with $d>0$ and $1\leq p\leq 2k-2$. 
Then 
\[v\sim v+2dn \cdot{\mathbf 1}\] 
for any $n\in\Z$. 
\end{lemma}

\begin{proof}
It suffices to consider the case $n=1$. 
Assume that $p$ is even;  
the case that $p$ is odd is treated analogously. 
Let $\beta=\sigma_{p+1}\cdots\sigma_{2}\sigma_{1}^{2}\sigma_{2}\cdots\sigma_p
\sigma_{p+1}^{-1}$. 
Figure~\ref{2d-reduction} shows that 
\[v\sim 
v\cdot\beta=(\underbrace{x+2d,\dots,x+2d}_{p},
x+\lambda d, 
\underbrace{x+d,\dots,x+d}_{2k-p-1}).
\]

\begin{figure}[htbp]
  \vspace{0.5em}
  \centering
    \begin{overpic}[]{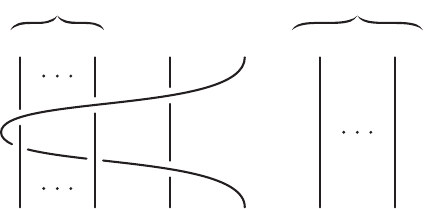}
      \put(26,97){{\footnotesize $p$}}
      \put(153,97){{\footnotesize $2k-p-1$}}
      \put(7,78){$x$}
      \put(43,78){$x$}
      \put(70,78){$x+\lambda d$}
      \put(106,78){$x+d$}
      \put(142,78){$x+d$}
      \put(178.5,78){$x+d$}
      \put(-2.5,-12){$x+2d$}
      \put(34,-12){$x+2d$}
      \put(70,-12){$x+\lambda d$}
      \put(106.5,-12){$x+d$}
      \put(142,-12){$x+d$}
      \put(178,-12){$x+d$}
    \end{overpic}
  \vspace{1em}
  \caption{The braid $\beta=\sigma_{p+1}\cdots\sigma_{2}\sigma_{1}^{2}\sigma_{2}
  \cdots\sigma_p\sigma_{p+1}^{-1}$}
  \label{2d-reduction}
\end{figure}
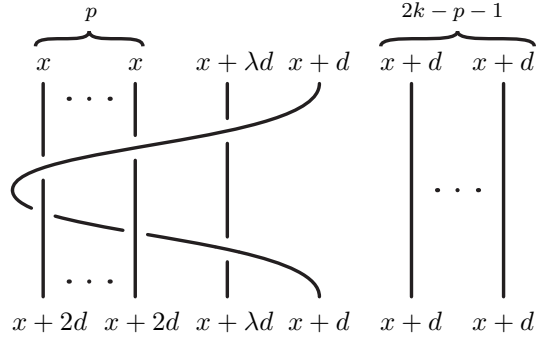

Now, we define 
$\widetilde{v},\widetilde{w}\in\Z^{2k-p+1}$ by 
\[\begin{split} 
&\widetilde{v} = 
(x+2d, x+\lambda d, 
\underbrace{x+d,\dots,x+d}_{2k-p-1}), \\
&\widetilde{w}=
(x+2d, x+(\lambda+2)d, \underbrace{x+3d,\dots,x+3d}_{2k-p-1}).  
\end{split}\]
A direct computation shows that  
\[\begin{split}
& \Delta(\widetilde{v})=\Delta(\widetilde{w})=
x+(3-\lambda)d,\\
& d(\widetilde{v})=d(\widetilde{w})=d, \\
& M(\widetilde{v})=M(\widetilde{w})=
\{x,x+\lambda d, 
\underbrace{x+d,\dots,x+d}_{2k-p-1}\}.
\end{split}
\]
Therefore, it follows from Theorem~\ref{thm47} that 
$\widetilde{v}\sim\widetilde{w}$, and hence, 
\[
v\sim(\underbrace{x+2d,\dots,x+2d}_{p-1},\widetilde{v})\sim
(\underbrace{x+2d,\dots,x+2d}_{p-1},\widetilde{w})
=v+2d\cdot{\mathbf 1}.\]
\end{proof}

\begin{proposition}\label{prop52}
Let $v$ be a nontrivial element in $\Z^{2k}$ $(k\geq 2)$, 
and let 
\[d=d(v)\mbox{ and }M(v)=\{\underbrace{r,\dots,r}_{p},
\underbrace{r+d,\dots,r+d}_{2k-p}\}\] 
for some $0\leq r<d$ and $1\leq p\leq 2k-1$. 
\begin{enumerate}
\item
If $2\leq p\leq 2k-1$, 
then there exists an even integer $\lambda$ such that 
\[v\sim(\underbrace{r,\dots,r}_{p-1}, 
r+\lambda d, 
\underbrace{r+d,\dots,r+d}_{2k-p}).\]
\item
If $1\leq p\leq 2k-2$, 
then there exists an odd integer $\lambda$ such that 
\[v\sim(\underbrace{r,\dots,r}_{p}, 
r+\lambda d, 
\underbrace{r+d,\dots,r+d}_{2k-p-1}).\]
\end{enumerate}
\end{proposition} 

\begin{proof}
Let $v=(a_1,\dots,a_{2k})$. 

(i) 
Since $p\geq 2$, Lemma~\ref{lem23} implies that 
there exist integers $1\leq i<j\leq 2k$ 
such that 
\[a_i\equiv a_j\equiv r \pmod{2d}.\] 
Therefore, by applying a Hurwitz action to $v$ if necessary, 
we may assume that 
\[a_{2k-1}\equiv a_{2k}\equiv r \pmod{2d}.\]
By Lemma~\ref{lem42}, we have  
\[(a_{2k-2},a_{2k-1}, a_{2k})\sim(x,y,y)\]
for some $x\leq y$ and $y\equiv r\pmod{2d}$. 
Moreover, it follows from Proposition~\ref{prop46}(i) that
\[\begin{split}
v& \sim 
(a_1,\dots,a_{2k-3},x,y,y)\\
& \sim
(\underbrace{\delta,\dots,\delta}_{p'},
\underbrace{\delta+d,\dots,\delta+d}_{2k-p'-1}, y)
\sim
(\underbrace{\delta,\dots,\delta}_{p'},y,
\underbrace{\delta+d,\dots,\delta+d}_{2k-p'-1})
\end{split}
\] 
for some $\delta\in\Z$ and $p'$ odd, 
where we use the fact that $2k-p'-1$ is even. 
Therefore, we have 
\[M(v)=\{\underbrace{\delta,\dots,\delta}_{p'},r,
\underbrace{\delta+d,\dots,\delta+d}_{2k-p'-1}\},\]
and hence, $\delta\equiv r$ or $r+d$ (mod~$2d$) 
by assumption. 

If $\delta\equiv r\pmod{2d}$, 
then $p'=p-1$ and $p$ is even. 
By Lemma~\ref{lem51}, 
we have 
\[v\sim
(\underbrace{\delta,\dots,\delta}_{p-1},y,
\underbrace{\delta+d,\dots,\delta+d}_{2k-p})
\sim
(\underbrace{r,\dots,r}_{p-1},y+2dn,
\underbrace{r+d,\dots,r+d}_{2k-p}), 
\]
where $n=(r-\delta)/2d\in{\Z}$. 
Since $y+2dn\equiv r\pmod{2d}$, the conclusion follows. 

If $\delta\equiv r+d\pmod{2d}$, 
then $p'=2k-p$ and $p$ is odd. 
By Lemma~\ref{lem51}, we have 
\begin{align*}
v&\sim 
(\underbrace{\delta,\dots,\delta}_{2k-p},y,
\underbrace{\delta+d,\dots,\delta+d}_{p-1})\\
&\sim
(\underbrace{r+d,\dots,r+d}_{2k-p},y+2dn, 
\underbrace{r+2d,\dots,r+2d}_{p-1}),
\end{align*}
where $n=(r+d-\delta)/2d\in{\Z}$. 
Figure~\ref{pf-normalform-even} shows that 
\[v\sim(\underbrace{r,\dots,r}_{2k-p},(2r-y)-2d(n-1),
\underbrace{r+d,\dots,r+d}_{2k-p}).\]
Since $(2r-y)-2d(n-1)\equiv 2r-r\equiv r\pmod{2d}$, 
the conclusion follows.

\begin{figure}[htbp]
  \vspace{1em}
  \centering
    \begin{overpic}[]{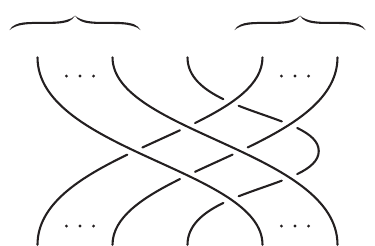}
      \put(24,115){{\footnotesize $2k-p$}}
      \put(135,115){{\footnotesize $p-1$}}
      \put(7,96){$r+d$}
      \put(43,96){$r+d$}
      \put(74.5,96){$y+2dn$}
      \put(113,96){$r+2d$}
      \put(149,96){$r+2d$}
      \put(158,44){$y+2dn$}
      \put(14,-9){$r$}
      \put(52,-9){$r$}
      \put(88,-12){$\uparrow$}
      \put(45,-25){$(2r-y)-2d(n-1)$}
      \put(115.5,-9){$r+d$}
      \put(152,-9){$r+d$}
    \end{overpic}
  \vspace{2em}
  \caption{Proof of Proposition~\ref{prop52}(i)}
  \label{pf-normalform-even}
\end{figure}

(ii) Since $2k-p\geq 2$, 
we may assume 
\[a_{2k-1}\equiv a_{2k}\equiv r+d \pmod{2d}\]
by Lemma~\ref{lem23}. 
It follows from Lemma~\ref{lem42} that 
\[(a_{2k-2},a_{2k-1}, a_{2k})\sim(x,y,y)\]
for some $x\leq y$ and $y\equiv r+d\pmod{2d}$. 
Moreover, it follows from Proposition~\ref{prop46}(i) that 
\[\begin{split}
v& \sim 
(a_1,\dots,a_{2k-3},x,y,y)\\
& \sim
(\underbrace{\delta,\dots,\delta}_{p'},
\underbrace{\delta+d,\dots,\delta+d}_{2k-p'-1}, y)
\sim
(\underbrace{\delta,\dots,\delta}_{p'},y,
\underbrace{\delta+d,\dots,\delta+d}_{2k-p'-1})
\end{split}
\] 
for some $\delta\in\Z$ and $p'$ odd, 
where we use the fact that 
$2k-p'-1$ is even. 
Therefore, we have 
\[M(v)=\{\underbrace{\delta,\dots,\delta}_{p'},r+d,
\underbrace{\delta+d,\dots,\delta+d}_{2k-p'-1}\},\]
and hence, $\delta\equiv r$ or $r+d\pmod{2d}$ 
by assumption.

If $\delta\equiv r\pmod{2d}$, 
then $p'=p$ and $p$ is odd. 
Therefore, 
the conclusion follows from Lemma~\ref{lem51}. 

If $\delta\equiv r+d\pmod{2d}$, 
then $p'=2k-p-1$ and $p$ is even. 
Similarly to the proof of (i), 
it follows that 
\[v\sim 
(\underbrace{r,\dots,r}_{2k-p-1},y',
\underbrace{r+d,\dots,r+d}_{p})
\]
for some $y'$ with $y'\equiv r+d\pmod{2d}$. 
\end{proof}

\begin{theorem}\label{thm53}
For $v,w\in\Z^{2k}$, 
the following are equivalent. 
\begin{enumerate}
\item 
$v\sim w$. 
\item 
$\Delta(v)=\Delta(w)$, 
$d(v)=d(w)$, and 
$M(v)=M(w)$. 
\end{enumerate}
\end{theorem}

\begin{proof}
\underline{(i)$\Rightarrow$(ii).} 
This follows from Lemma~\ref{lem22}. 

\underline{(ii)$\Rightarrow$(i).} 
If $v$ is trivial, 
then since $d(w)=d(v)=0$, so is $w$. 
Therefore, since $M(v)=M(w)$, it follows that $v=w$. 

Assume that $v$ and $w$ are nontrivial. 
Let $d=d(v)=d(w)>0$. 
Since $M(v)=M(w)$, it follows from Proposition~\ref{prop52} that 
\begin{align*}
&v\sim(\underbrace{r,\dots,r}_{p}, r+\lambda d,
\underbrace{r+d,\dots,r+d}_{2k-p-1}), \\
&w\sim(\underbrace{r,\dots,r}_{p}, r+\lambda' d, 
\underbrace{r+d,\dots,r+d}_{2k-p-1})
\end{align*}
for some $r,p,\lambda,\lambda'\in\Z$. 
Moreover, since $\Delta(v)=\Delta(w)$, 
we have $\lambda=\lambda'$, 
and hence, $v\sim w$. 
\end{proof}

\begin{proposition}\label{prop54}
A complete system of representatives 
of $\Z^{2k}/\sim$ is given by 
\[\{(x,\dots,x)\in{\Z}^{2k}\mid x\in{\Z}\} 
\sqcup D\sqcup\bigsqcup_{p=1}^{2k-2}E_p,\]
where $D$ and $E_p$ $(1\leq p\leq 2k-2)$ 
are the subsets of ${\Z}^{2k}$ defined by 
\begin{align*}
D&=\left\{\left.
(x, x+(2n-1)(y-x),
\underbrace{y,\dots,y}_{2k-2})
\in{\Z}^{2k}
\right|\ 
0\leq 2x<y, \ n\in{\Z}\right\}, \medskip\\
E_p&=\left\{\left.
(\underbrace{x,\dots,x}_{p}, x+2n(y-x),
\underbrace{y,\dots,y}_{2k-p-1})
\in{\Z}^{2k}
\right| \ 
0\leq 2x<y, \ n\in{\Z}
\right\}.
\end{align*}
\end{proposition}

\begin{proof}
Let $v$ be a nontrivial element of ${\Z}^{2k}$ 
with 
\[M(v)=\{\underbrace{r,\dots,r}_{p},
\underbrace{r+d,\dots,r+d}_{2k-p}\}\] 
for some $r,d,p\in\Z$ with $0\leq r<d$ and $1\leq p\leq 2k-1$. 
By Proposition~\ref{prop52}, 
$v$ is equivalent to 
\begin{enumerate}
\item
$(r,r+\lambda d,\underbrace{r+d,\dots,r+d}_{2k-2})$ 
for some $\lambda$ odd if $p=1$, 
\item
$(\underbrace{r,\dots,r}_{p-1},r+\lambda d,\underbrace{r+d,\dots,r+d}_{2k-p})$ 
for some $\lambda$ even if $2\leq p\leq 2k-2$. 
\end{enumerate}
Since $0\leq 2r<r+d$, the element in (i) (resp. (ii)) belongs to $D$ 
(resp. $E_p$). 

Let $v_1$ and $v_{2}$ be nontrivial elements 
in $D\sqcup\bigsqcup_{p=1}^{2k-2}E_p$ such that 
\[v_i=(\underbrace{x_i,\dots,x_i}_{p_i}, x_i+\ell_i(y_i-x_i),
\underbrace{y_i,\dots,y_i}_{2k-p_i-1}) \quad 
(i=1,2).\] 
Here, $x_i$ and $y_i$ satisfy $0\leq 2x_i<y_i$, 
and $\ell_i$ is odd (resp. even) if $v_i\in D$ 
(resp. $v_i\in\bigsqcup_{p=1}^{2k-2}E_p$). 
A direct computation shows that 
\begin{align*}
&\Delta(v_i)=
\begin{cases}
-\ell_i(y_i-x_i) & \mbox{for $p_i$ odd}, \\
(\ell_i-1)(y_i-x_i) & \mbox{for $p_i$ even},
\end{cases}\\
&d(v_i)=y_i-x_i, \\
&M(v_i)=
\begin{cases}
\{x_i,\underbrace{y_i,\dots,y_i}_{2k-1}\} & \mbox{for $v_i\in D$}, \\
\{\underbrace{x_i,\dots,x_i}_{p_i+1},
\underbrace{y_i,\dots,y_i}_{2k-p_i-1}\} & 
\mbox{for $v_i\in \bigsqcup_{p=1}^{2k-2}E_p$}.
\end{cases}
\end{align*}
We remark that 
\[
0\leq x_{i}<d(v_{i}) \text{ and } 
d(v_{i})\leq y_{i}<2d(v_{i}). 
\]
Assume that $v_1\sim v_2$. 
Since $d(v_1)=d(v_2)$ and $M(v_1)=M(v_2)$, 
we have $x_1=x_2$, $y_1=y_2$, and $p_1=p_2$. 
In particular, both $v_1$ and $v_2$ lie in the same subset 
among $D,E_1,\dots,E_{2k-2}$. 
Moreover, since $\Delta(v_1)=\Delta(v_2)$, 
we have $\ell_1=\ell_2$, 
and hence, $v_1=v_2$. 
\end{proof}

By Proposition~\ref{prop49}, any element in ${\Z}^{2k-1}$ 
is equivalent to one of the form $(x,\dots,x,y,\dots,y)$. 
Although this property does not hold in general for $m=2k$, 
we have the following.  

\begin{proposition}\label{prop55}
For a nontrivial element $v\in\Z^{2k}$ $(k\geq2)$, 
the following are equivalent. 
\begin{enumerate}
\item 
$v\sim(x,\dots,x,y,\dots,y)$ for some $x\ne y$. 
\item 
$\Delta(v)\in\{0,\pm d(v)\}$. 
\end{enumerate}
In particular, 
$v\sim (\underbrace{x,\dots,x}_{\rm even}, 
\underbrace{y,\dots,y}_{\rm even})$ 
if and only if $\Delta(v)=0$. 
\end{proposition}

\begin{proof}
\underline{(i)$\Rightarrow$(ii).} 
Assume that 
\[v\sim(\underbrace{x,\dots,x}_p,\underbrace{y,\dots,y}_{2k-p})\]
for some $1\leq p\leq 2k-1$. 
By Lemma~\ref{lem22}(i) and (ii), we have 
\begin{align*}
&\Delta(v)=
\begin{cases}
x-y & \mbox{for $p$ odd}, \\
0 & \mbox{for $p$ even}, 
\end{cases}\\
&d(v)=|y-x|, 
\end{align*}
and hence, the conclusion (ii) follows. 

\underline{(ii)$\Rightarrow$(i).} 
By Proposition~\ref{prop54}, 
we may assume that 
\[v=(\underbrace{x,\dots,x}_p, x+\ell(y-x), 
\underbrace{y,\dots,y}_{2k-p-1})\] 
for some $x,y,\ell,p\in\Z$ with 
$0\leq 2x<y$ and $1\leq p\leq 2k-2$.

Firstly, we consider the case that $p$ is even. 
Since $\Delta(v)=(\ell-1)(y-x)$, 
we have $\ell=0,1,2$ by assumption. 
For $\ell=0$ or $1$, the conclusion (i) already follows. 
For $\ell=2$, 
since we have 
\[\Delta(v)=d(v)=y-x \mbox{ and }
M(v)=\{\underbrace{x,\dots,x,x}_{p+1},
\underbrace{y,\dots,y}_{2k-p-1}\},\] 
it follows from Theorem~\ref{thm53} that 
\[v= 
(\underbrace{x,\dots,x}_{p}, 
x+2(y-x), 
\underbrace{y,\dots,y}_{2k-p-1})
\sim (\underbrace{y,\dots,y}_{2k-p-1},
\underbrace{x,\dots,x,x}_{p+1}).
\]

Secondly, we consider the case that $p$ is odd. 
Since $\Delta(v)=-\ell(y-x)$,  
we have $\ell=0,\pm 1$ by assumption. 
For $\ell=0$ or $1$, 
the conclusion (i) already follows. 
For $\ell=-1$, 
since we have 
\[\Delta(v)=d(v)=y-x \mbox{ and }
M(v)=\{\underbrace{x,\dots,x}_{p}, 
\underbrace{y,y,\dots,y}_{2k-p}\},\] 
it follows from Theorem~\ref{thm53} that 
\begin{align*}
v&= 
(\underbrace{x,\dots,x}_p, 
x-(y-x), 
\underbrace{y,\dots,y}_{2k-p-1})\\
&=
(\underbrace{x,\dots,x}_p, 
y-2(y-x), 
\underbrace{y,\dots,y}_{2k-p-1})
\sim (\underbrace{y,y,\dots,y}_{2k-p},
\underbrace{x,\dots,x}_{p}).
\end{align*}

In particular, 
$\Delta(v)=0$ if and only if 
\begin{itemize}
\item
$p$ is even and $\ell=1$, and hence,  
$v=(\underbrace{x,\dots,x}_p,\underbrace{y,\dots,y}_{2k-p})$, or 
\item
$p$ is odd and $\ell=0$, and hence, 
$v=(\underbrace{x,\dots,x}_{p+1},\underbrace{y,\dots,y}_{2k-p-1})$. 
\end{itemize}
In both cases, 
the numbers of $x$'s and $y$'s appearing in $v$ are even. 
\end{proof}

\section{Pure braid groups}\label{sec6}

A {\it pure $m$-braid} is an $m$-braid $\beta\in B_m$ 
such that the associated permutation 
$\pi_\beta\in S_m$ is the identity $e$; 
that is, for each $i=1,\dots,m$, 
the $i$th top point of $\beta$ from the left 
connects to the $i$th bottom point 
via a string.
Let $P_m$ denote the pure $m$-braid group, 
which is the normal subgroup of $B_{m}$ 
consisting of all pure $m$-braids. 
We denote by $\stackrel{P}{\sim}$ 
the equivalence relation on $\Z^{m}$ induced by 
the Hurwitz action of $P_m\subset B_m$. 
By definition, $v\stackrel{P}{\sim}w$ 
implies $v\sim w$. 

For $v=(a_1,\dots,a_m)\in{\Z}^m$, 
let $M^*(v)$ denote 
the ordered $m$-tuple 
consisting of the congruence classes of 
$a_1,\dots,a_m$ modulo $2d(v)$; 
that is, 
\[M^*(v)=(a_1,\dots,a_m)\in({\Z}/2d(v){\Z})^m.\]
By definition, $v=(a_1,\dots,a_m)$ and 
$w=(b_1,\dots,b_m)$ satisfy that $M^*(v)=M^*(w)$ 
if and only if $a_i\equiv b_i\pmod{2d(v)}$ 
for any $i=1,\dots,m$. 

\begin{lemma}\label{lem61}
If $v\stackrel{P}{\sim}w$, then $M^*(v)=M^*(w)$. 
\end{lemma}

\begin{proof}
This is a direct consequence of Lemma~\ref{lem23}. 
\end{proof}

For integers $p$ and $q$ with $1\leq p<q\leq m$, 
let $S_{m}(p,q)$ denote the subgroup of $S_m$ 
consisting of all permutations 
on $\{p,\dots,q\}$;  
that is, 
\[S_m(p,q)=\{\pi\in S_m\mid 
\pi(i)=i \mbox{ for }1\leq i<p\mbox{ and }
q<i\leq m\}.\]

\begin{lemma}\label{lem62} 
Let $v\in\Z^{m}$ $(m\geq3)$ be an element of the form 
\[v=(\underbrace{x,\dots,x}_p,x+\ell(y-x), 
\underbrace{y,\dots,y}_{m-p-1}).\]
for some $x,y,\ell,p\in\Z$ with 
$0\leq 2x<y$ and $1\leq p\leq m-2$. 
\begin{itemize}
\item[{\rm (i)}] 
If $\ell$ is even, then for any 
$\pi\in S_m(1,p+1)$ and 
$\pi'\in S_m(p+2,m)$, 
there exists an $m$-braid $\beta\in B_{m}$ 
such that \[v\cdot \beta=v\mbox{ and }\pi_{\beta}=\pi\pi'.\]
\item[{\rm (ii)}] 
If $\ell$ is odd, then for any 
$\pi\in S_m(1,p)$ and 
$\pi'\in S_m(p+1,m)$, 
there exists an $m$-braid $\beta\in B_{m}$ 
such that \[v\cdot \beta=v\mbox{ and }\pi_{\beta}=\pi\pi'.\]
\end{itemize}
\end{lemma}

\begin{proof}
We only prove (i); 
the proof for (ii) follows analogously. 
It suffices to consider the case that  
$\pi$ is an adjacent transposition $(i\ i+1)$ 
for $1\leq i\leq p$ 
and $\pi'=e$. 
For $1\leq i\leq p-1$, 
the conclusion follows 
by taking $\beta=\sigma_i$. 

We consider the case $i=p$; that is, $\pi=(p\ p+1)$ and $\pi'=e$. 
Since $\ell$ is even, it follows from Theorem~\ref{thm47} that 
\[(x,x+\ell(y-x),y)\sim
(x-\ell(y-x),x-\ell(y-x),y-\ell(y-x))\in{\Z}^3.\] 
Therefore, there exists an $m$-braid $\gamma\in B_{m}$ 
such that 
\[v\cdot \gamma=
(\underbrace{x,\dots,x}_{p-1},x-\ell(y-x),x-\ell(y-x),y-\ell(y-x),
\underbrace{y,\dots,y}_{m-p-2})\]
and $\pi_{\gamma}=e\mbox{ or }(p\ p+1)$ by Lemma~\ref{lem23}. 
Let $\beta=\gamma\sigma_p\gamma^{-1}$. 
Then it follows that 
\begin{align*}
&v\cdot\beta=((v\cdot\gamma)\cdot\sigma_p)\cdot\gamma^{-1}
=(v\cdot\gamma)\cdot\gamma^{-1}=v, \\
&\pi_{\beta}=\pi_{\gamma}(p\ p+1)\pi_{\gamma}^{-1}=(p\ p+1)=\pi\pi'.
\end{align*}
\end{proof}

\begin{figure}[htbp]
  \vspace{1em}
  \centering
    \begin{overpic}[]{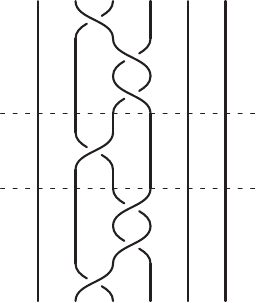}
      \put(4,115){$\gamma$}
      \put(16,149){$0$}
      \put(33.5,149){$0$}
      \put(52,149){$2$}
      \put(70,149){$1$}
      \put(88,149){$1$}
      \put(106,149){$1$}
      \put(3,70){$\sigma_{2}$}
      \put(21,82){$-2$}
      \put(56.5,82){$-2$}
      \put(75,82){$-1$}
      \put(21,58){$-2$}
      \put(56.5,58){$-2$}
      \put(75,58){$-1$}
      %
      \put(-2,25){$\gamma^{-1}$}
      \put(16,-10){$0$}
      \put(33.5,-10){$0$}
      \put(52,-10){$2$}
      \put(70,-10){$1$}
      \put(88,-10){$1$}
      \put(106,-10){$1$}
      \end{overpic}
  \vspace{1em}
  \caption{The $6$-braid 
  $\beta=(\sigma_2^{-1}\sigma_3^{-2})\sigma_2(\sigma_2^{-1}\sigma_3^{-2})^{-1}$
  with $\pi_{\beta}=(2\ 3)$}
  \label{ex-lem-permutation}
\end{figure}
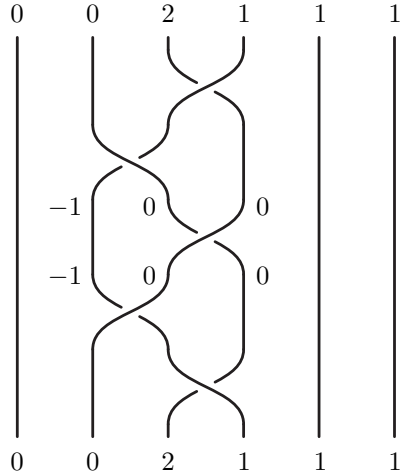

For example, 
we consider the element 
$v=(0,0,2,1,1,1)\in\Z^6$ 
as in Lemma~\ref{lem61} for 
$x=0$, $y=1$, $\ell=2$, 
and $p=2$. 
For $\pi=(2\ 3)\in S_6(1,3)$ and $\pi'=e\in S_6(4,6)$, 
the $6$-braid 
$\gamma=\sigma_2^{-1}\sigma_3^{-2}$ 
satisfies that 
\[v\cdot\gamma=(0,-2,-2,-1,1,1) \mbox{ and }\pi_{\gamma}=(2\ 3).\]
Therefore, 
$\beta=\gamma\sigma_2\gamma^{-1}$ shows that  
$v\cdot\beta=v$ and $\pi_{\beta}=(2\ 3)=\pi\pi'$. 
See Figure~\ref{ex-lem-permutation}. 

\begin{theorem}\label{thm63}
For $v,w\in\Z^m$ $(m\geq2)$, 
the following are equivalent. 
\begin{enumerate}
\item 
$v\stackrel{P}{\sim} w$. 
\item 
$\Delta(v)=\Delta(w)$, 
$d(v)=d(w)$, and 
$M^*(v)=M^*(w)$. 
\end{enumerate}
\end{theorem}

\begin{proof}
\underline{(i)$\Rightarrow$(ii).}
This follows from Lemmas~\ref{lem22}(i), (ii), and \ref{lem61}. 

\underline{(ii)$\Rightarrow$(i).} 
If $v$ is trivial, 
then since $d(w)=d(v)=0$, so is $w$. 
Therefore, since $M^*(v)=M^*(w)$, we have $v=w$. 

Assume that $v$ and $w$ are nontrivial. 
By assumption, we have $M(v)=M(w)$, 
and hence, $v\sim w$ by Theorem~\ref{thm11}. 

\underline{$m=2$.}
Let $v=(a_1,a_2)$ and $w=(b_1,b_2)$. 
Note that $d(v)=|a_2-a_1|=|\Delta(v)|$. 
Since $\Delta(v)=\Delta(w)$ and $M^*(v)=M^*(w)$, 
we have 
\[a_1-a_2=b_1-b_2\mbox{ and }a_1-b_1=2n\Delta(v)\] 
for some $n\in\Z$. 
Then we have $a_2-b_2=2n\Delta(v)$, and hence, 
\[w=(a_1-2n\Delta(v),a_2-2n\Delta(v))
=v\cdot\sigma_1^{2n}.\]
Since $\sigma_1^{2n}$ is a pure $2$-braid, 
the conclusion $v\stackrel{P}{\sim}w$ follows. 

\underline{$m=2k-1$ $(k\geq 2)$.}
By Proposition~\ref{prop49}, there exist 
$\beta_1,\beta_2\in B_{2k-1}$ such that 
\[v\cdot \beta_1=(\underbrace{x,\dots,x}_{2p-1}, 
\underbrace{y,\dots,y}_{2k-2p}) 
=w\cdot\beta_2.\]
Since $M^*(v)=M^*(w)$, 
it follows from Lemma~\ref{lem23} that 
there exist $\pi\in S_{2k-1}(1,2p-1)$ and $\pi'\in S_{2k-1}(2p,2k-1)$ 
such that $\pi\pi'=\pi_{\beta_1}^{-1}\pi_{\beta_2}$. 
Hence, by Lemma~\ref{lem62}, 
there exists $\beta_3\in B_{2k-1}$ such that 
\[(\underbrace{x,\dots,x}_{2p-1}, 
\underbrace{y,\dots,y}_{2k-2p}) \cdot \beta_3=
(\underbrace{x,\dots,x}_{2p-1}, 
\underbrace{y,\dots,y}_{2k-2p})\] 
and $\pi_{\beta_3}=\pi\pi'$. 
Then we have 
$v\cdot(\beta_1\beta_3\beta_2^{-1})=w$ and 
\[\pi_{\beta_1\beta_3\beta_2^{-1}}
=\pi_{\beta_1}\pi_{\beta_3}\pi_{\beta_2}^{-1}
=\pi_{\beta_1}(\pi_{\beta_1}^{-1}
\pi_{\beta_2})\pi_{\beta_2}^{-1}=e.\]
Since 
$\beta_1\beta_3\beta_2^{-1}$ is a pure $(2k-1)$-braid,  
the conclusion $v\stackrel{P}{\sim}w$ follows.

\underline{$m=2k$ $(k\geq 2)$.} 
By Proposition~\ref{prop54}, 
there exist $\beta_1,\beta_2\in B_{2k}$ such that 
\[
v\cdot\beta_1=
(\underbrace{x,\dots,x}_p,x+\ell(y-x),
\underbrace{y,\dots,y}_{2k-p-1})
=w\cdot\beta_2.\]
We consider only the case that $\ell$ is even; 
the case that $p=1$ and $\ell$ is odd is treated analogously. 
Since $M^*(v)=M^*(w)$, it follows from Lemma~\ref{lem23} that 
there exist $\pi\in S_{2k}(1,p+1)$ and $\pi'\in S_{2k}(p+2,2k)$ 
such that $\pi\pi'=\pi_{\beta_1}^{-1}\pi_{\beta_2}$. 
Hence, by Lemma~\ref{lem62}, 
there exists $\beta_3\in B_{2k}$ such that 
\[(\underbrace{x,\dots,x}_p,x+\ell(y-x),
\underbrace{y,\dots,y}_{2k-p-1})\cdot\beta_3
=(\underbrace{x,\dots,x}_p,x+\ell(y-x),
\underbrace{y,\dots,y}_{2k-p-1})\] 
and $\pi_{\beta_3}=\pi\pi'$. 
Then we have 
$v\cdot(\beta_1\beta_3\beta_2^{-1})=w$ and 
$\pi_{\beta_1\beta_3\beta_2^{-1}}=e$. 
Since $\beta_1\beta_3\beta_2^{-1}$ is a pure $2k$-braid, 
the conclusion $v\stackrel{P}{\sim}w$ follows. 
\end{proof}

For example, three elements 
\[v=(1,-5,4), \ w=(10,7,7),\mbox{ and }
u=(7,7,10)\in\Z^3\]
satisfy that 
\[\begin{split}
& \Delta(v)=\Delta(w)=\Delta(u)=10, \ 
d(v)=d(w)=d(u)=3, \\ 
&M(v)=M(w)=M(u)=\{1,1,4\},  \\
&M^*(v)=M^*(u)
=(1,1,4)\neq(4,1,1)
=M^*(w). 
\end{split}\] 
Therefore, it follows from Theorems~\ref{thm47} and \ref{thm63} that 
\[v\sim w, \ v\stackrel{P}{\not\sim}w, 
\text{ and }v\stackrel{P}{\sim}u.\]

\section{Virtual braid groups}\label{sec7}

\begin{figure}[htbp]
  \vspace{1em}
  \begin{center}
   \hspace{1em}
    \begin{overpic}[]{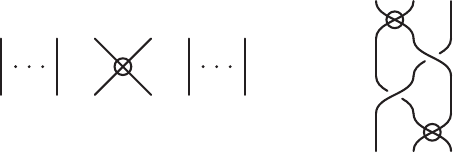}
      \put(-25,39){$\tau_{i}=$}
      \put(-3.5,60){$a_{1}$}
      \put(15,60){$a_{i-1}$}
      \put(39,60){$a_{i}$}
      \put(62,60){$a_{i+1}$}
      \put(84,60){$a_{i+2}$}
      \put(114,60){$a_{m}$}
      \put(-3.5,18){$a_{1}$}
      \put(15,18){$a_{i-1}$}
      \put(39,18){$a_{i+1}$}
      \put(70,18){$a_{i}$}
      \put(84,18){$a_{i+2}$}
      \put(114,18){$a_{m}$}
      \put(179,77){$1$}
      \put(193,77){$-5$}
      \put(214,77){$4$}
      \put(174,-11){$-2$}
      \put(197,-11){$1$}
      \put(215,-11){$1$}
    \end{overpic}
  \end{center}
  \vspace{1em}
  \caption{The Hurwitz action of $VB_{m}$ on $\Z^{m}$}
  \label{Hurwitz-action-virtual}
\end{figure}

Let $VB_m$ $(m\geq 2)$ 
denote the {\it virtual $m$-braid group} 
with the standard generators 
\[\sigma_1,\dots,\sigma_{m-1} \text{ and }
\tau_1,\dots,\tau_{m-1},\] 
where $\sigma_i$ represents a classical crossing 
between the $i$th and $(i+1)$st strings 
(see again the left of Figure~\ref{Hurwitz-action}), 
and $\tau_i$ represents a virtual crossing 
between the $i$th and $(i+1)$st strings
for $i=1,\dots,m-1$. 
See~\cite{Kam03,KL} for more details of virtual braid groups. 
The set $\Z^m$ admits the Hurwitz action of $VB_m$, 
in which the generators 
$\sigma_{i}$ act as defined in Section~\ref{sec2}, 
while the generators $\tau_{i}$ act by
\[(a_{1},\dots,a_{i},a_{i+1},\dots,a_{m})\cdot\tau_i=(a_1,\dots,a_{i-1},a_{i+1},a_i,a_{i+2},\dots,a_m).\]
See the left of Figure~\ref{Hurwitz-action-virtual}. 
We denote by $v\stackrel{V}{\sim}w$ 
if $v\cdot \beta=w$ holds for some $\beta\in VB_m$. 
The right of the figure shows that 
\[
(1,-5,4)\cdot\tau_{1}\sigma_{2}\sigma_{1}^{-1}\tau_{2}=(-2,1,1)\in\Z^{3}, 
\]
which means that $(1,-5,4)\stackrel{V}{\sim}(-2,1,1)$. 

For a virtual $m$-braid $\beta\in VB_m$, 
we denote by $\pi_{\beta}\in S_m$ 
the permutation associated with $\beta$. 
The {\it virtual pure $m$-braid group} is defined by 
\[VP_m=\{\beta\in VB_m\mid 
\pi_{\beta}=e\},\]
which is a normal subgroup of $VB_m$. 
As with $VB_m$, the group $VP_m$ acts on $\Z^m$. 
For $v,w\in\Z^m$, 
we write $v\stackrel{VP}{\sim}w$ 
if there exists $\beta\in VP_m$ 
such that $v\cdot\beta=w$. 
By definition, $v\stackrel{VP}{\sim}w$ implies  
$v\stackrel{V}{\sim}w$.

The following two lemmas can be proved 
in the same way as Lemmas~\ref{lem22} and~\ref{lem23}, 
and we omit the proofs. 

\begin{lemma}\label{lem71}
If $v\stackrel{V}{\sim}w\in\Z^{m}$ $(m\geq2)$, 
then $d(v)=d(w)$ and $M(v)=M(w)$.
\qed
\end{lemma}

\begin{lemma}\label{lem72}
Suppose that 
$v=(a_1,\dots,a_m)$,  
$w=(b_1,\dots,b_m)\in\Z^m$ $(m\geq2)$, 
and $\beta\in VB_m$ satisfy $v\cdot \beta=w$. 
Then we have 
\[b_{\pi_\beta(k)}\equiv a_k\pmod{2d(v)}\]
for any $k=1,\dots,m$. 
Therefore, if $v\stackrel{VP}{\sim}w$, 
then 
$M^*(v)=M^*(w)$. 
\qed
\end{lemma}

\begin{lemma}\label{lem73}
If two elements $v$ and $w\in\Z^{m}$ $(m\geq2)$
satisfy that $d(v)=d(w)$ and $M(v)=M(w)$, 
then $\Delta(v)\equiv\Delta(w)\pmod{2d(v)}$. 
\end{lemma}

\begin{proof}
We may assume that $v$ and $w$ are nontrivial. 
By Lemma~\ref{lem21}(iv), we have 
\[M(v)=M(w)=\{\underbrace{r,\dots,r}_p,
\underbrace{r+d,\dots,r+d}_{m-p}\},\] 
where $d=d(v)=d(w)$,  
$0\leq r<d$, and $1\leq p\leq m-1$. 
Then it follows that 
\[\Delta(v)\equiv\Delta(w)\equiv
\left\{
\begin{tabular}{rl}
$r+(m-p)d \pmod{2d}$ & for $m$ odd, \\
$(m-p)d \pmod{2d}$ & for $m$ even.
\end{tabular}
\right.\]
\end{proof}

\begin{lemma}\label{lem74}
For any $x,y,\ell,n\in\Z$, 
it holds that 
\[v=(x,x+\ell(y-x), y)
\stackrel{VP}{\sim} 
(x,x+(\ell+2n)(y-x), y). 
\]
\end{lemma}

\begin{figure}[htbp]
  \centering
  \vspace{2.5em}
    \begin{overpic}[]{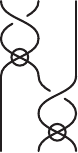}
      \put(-2.5,77){$x$}
      \put(16,79){$\downarrow$}
      \put(-6,92){$x+\ell(y-x)$}
      \put(34,78){$y$}
      \put(-58,52){$x-\ell(y-x)$}
      \put(-2.5,-8){$x$}
      \put(16,-10){$\uparrow$}
      \put(-20,-23){$x+(\ell+2)(y-x)$}
      \put(34,-8){$y$}
          \end{overpic}
  \vspace{2em}
  \caption{Proof of Lemma~\ref{lem74}}
  \label{2d-reduction-virtual}
\end{figure}

\begin{proof}
It suffices to consider the case $n=1$. 
Figure~\ref{2d-reduction-virtual} shows that 
\begin{align*}
(x,x+\ell(y-x),y)
&\stackrel{VP}{\sim}
(x,x+\ell(y-x),y)\cdot\sigma_1^{-1}\tau_1\sigma_2\tau_2\\
&=(x,x+(\ell+2)(y-x), y).
\end{align*}
\end{proof}

\begin{lemma}\label{lem75}
For a nontrivial element $v\in\Z^m$ $(m\geq3)$ 
and an integer $n\in\Z$, 
there exists an element $v'\in\Z^m$ such that 
\begin{enumerate}
\item
$v'\stackrel{VP}{\sim}v$ and 
\item
$\Delta(v')=\Delta(v)+2nd(v)$.
\end{enumerate}
\end{lemma}

\begin{proof}
By Propositions~\ref{prop49} and \ref{prop54}, 
there exists $\beta_1\in B_m$ such that  
\[v\sim 
v\cdot\beta_1=(\underbrace{x,\dots,x}_{p}, x+\ell(y-x), 
\underbrace{y,\dots,y}_{m-p-1})\] 
for some $x,y,\ell,p\in\Z$ with $x<y$ and $1\leq p\leq m-2$. 
By Lemma~\ref{lem74}, there exists $\beta_2\in VP_m$ 
such that 
\[v\cdot\beta_1\stackrel{VP}{\sim}(v\cdot\beta_1)\cdot\beta_2=
(\underbrace{x,\dots,x}_{p}, x+(\ell+2(-1)^pn)(y-x), 
\underbrace{y,\dots,y}_{m-p-1}).\] 
Then we have
\[\Delta(v\cdot\beta_1\beta_2)=\Delta(v)+2n(y-x).\]
We define $v'=v\cdot\beta_1\beta_2\beta_1^{-1}$. 
Since $\beta_1\beta_2\beta_1^{-1}\in VP_m$, 
we have $v\stackrel{VP}{\sim}v'$. 
Moreover, it follows from Lemma~\ref{lem22}(i) that 
\[\Delta(v')=\Delta(v\cdot\beta_1\beta_2)=\Delta(v)+2n(y-x).\]
Since $d(v)=y-x$, the conclusion follows. 
\end{proof}

\begin{theorem}\label{thm76}
For $v,w\in\Z^m$ $(m\geq2)$, 
we have the following. 
\begin{enumerate}
\item 
$v\stackrel{V}{\sim} w$ if and only if 
$d(v)=d(w)$ and 
$M(v)=M(w)$. 
\item
$v\stackrel{VP}{\sim}w$ if and only if 
$d(v)=d(w)$ and 
$M^*(v)=M^*(w)$. 
\end{enumerate}
\end{theorem}

\begin{proof}
(i) ($\Rightarrow$) 
This follows from Lemma~\ref{lem71}. 

($\Leftarrow$) 
If $v$ is trivial, 
then since $d(w)=d(v)=0$, so is $w$. 
Therefore, since $M(v)=M(w)$, it follows that $v=w$.  
Assume that $v$ and $w$ are nontrivial. 

\underline{$m=2$.}
Since $d(v)=d(w)$, we have $\Delta(w)=\pm\Delta(v)$. 
If $\Delta(w)=\Delta(v)$, 
then it follows from Theorem~\ref{thm11} that $v\sim w$. 
If $\Delta(w)=-\Delta(v)$, 
then $\Delta(v\cdot\tau_1)=-\Delta(v)=\Delta(w)$. 
Therefore, we have 
$v\stackrel{V}{\sim}v\cdot\tau_1\sim w$.

\underline{$m\geq 3$.}
By Lemma~\ref{lem73}, 
we may assume that $\Delta(w)=\Delta(v)+2nd(v)$ 
for some $n\in{\Z}$. 
By Lemma~\ref{lem75}, 
there exists an element $v'\in{\Z}^m$ such that 
\[v'\stackrel{VP}{\sim}v \mbox{ and }
\Delta(v')=\Delta(v)+2nd(v).\]
By assumption and Lemma~\ref{lem71}, we have 
\[\Delta(v')=\Delta(w), \ d(v')=d(v)=d(w), \mbox{ and }
M(v')=M(v)=M(w).\] 
Therefore, it follows from Theorem~\ref{thm11} that $v'\sim w$, 
and hence $v\stackrel{V}{\sim}w$. 

(ii) 
$(\Rightarrow$) 
This follows from Lemmas~\ref{lem71} and \ref{lem72}. 

$(\Leftarrow$) 
The proof proceed in the same manner as that of (i), 
with the condition $M(v)=M(w)$ replaced by $M^*(v)=M^*(w)$. 
By Theorem~\ref{thm63}, we have $v'\stackrel{P}{\sim}w$, 
and hence, $v\stackrel{VP}{\sim}w$. 
\end{proof}

\begin{remark}\label{rem77}
(i) For $\beta\in VB_m$, 
if there exists an element $v\in\Z^{m}$ with 
$\Delta(v)\neq\Delta(v\cdot\beta)$, 
then it follows from Lemma~\ref{lem22}(i) that $\beta\not\in B_m$. 
The virtual $3$-braid on the right of 
Figure~\ref{Hurwitz-action-virtual} is such an example. 

(ii) The \textit{welded $m$-braid group} $WB_{m}$~\cite{Kam03,KL}, 
also known as the \textit{braid-permutation group}~\cite{FRR}, 
is a quotient of the virtual $m$-braid group $VB_{m}$. 
The Hurwitz action of $WB_m$ is defined via 
the Hurwitz action of an arbitrary representative in $VB_{m}$; 
this action is well-defined. 
Consequently, the equivalence relation on $\Z^{m}$ induced by the action of $WB_m$ coincides with that induced by $VB_{m}$. 
\end{remark}

\begin{proposition}\label{prop78}
A complete system of representatives 
of $\Z^{m}/\stackrel{V}{\sim}$ $(m\geq2)$ is given by 
\[\{(x,\dots,x)\in{\Z}^{m}\mid 
x\in{\Z}\}\sqcup \bigsqcup_{p=1}^{m-1}F_p,\] 
where $F_p$ $(1\leq p\leq m-1)$ is the subset of ${\Z}^{m}$ defined by 
\[F_p=\left\{\left.
(\underbrace{x,\dots,x}_{p},
\underbrace{y,\dots,y}_{m-p})
\in{\Z}^{m}
\right|\ 0\leq 2x<y
\right\}.\] 
\end{proposition}

\begin{proof}
Let $v\in{\Z}^m$ be a nontrivial element, 
and let $d=d(v)$ and 
\[M(v)=\{\underbrace{r,\dots,r}_{p},
\underbrace{r+d,\dots,r+d}_{m-p}\}\]
for some $0\leq r<d$ and $1\leq p\leq m-1$. 
By Theorem~\ref{thm76}(i), 
we have 
\[v\stackrel{V}{\sim}
(\underbrace{r,\dots,r}_{p}, 
\underbrace{r+d,\dots,r+d}_{m-p}).\]
Since $0\leq 2r<r+d$, 
this element belongs to $F_p$.

Let $v_1$ and $v_{2}$ be nontrivial elements 
in $\bigsqcup_{p=1}^{m-1}F_p$ such that 
\[v_i=(\underbrace{x_i,\dots,x_i}_{p_i}, 
\underbrace{y_i,\dots,y_i}_{m-p_i}) \quad 
(i=1,2).\] 
Here, $x_i$ and $y_i$ satisfy $0\leq 2x_i<y_i$. 
A direct computation shows that 
\[d(v_i)=y_i-x_i \mbox{ and }
M(v_i)=
\{\underbrace{x_i,\dots,x_i}_{p_i},
\underbrace{y_i,\dots,y_i}_{m-p_i}\}.
\]
We remark that 
\[
0\leq x_{i}<d(v_{i}) \text{ and } 
d(v_{i})\leq y_{i}<2d(v_{i}). 
\]
Assume that $v_1\sim v_2$. 
Since $d(v_1)=d(v_2)$ and $M(v_1)=M(v_2)$, 
we have $x_1=x_2$, $y_1=y_2$, and $p_1=p_2$, 
and hence, $v_1=v_2$. 
\end{proof}

For a nontrivial element $v\in\Z^{m}$, 
the positive integer $d(v)$ can be uniquely written in the form 
\[
d(v)=2^{k(v)}(2n+1)
\] 
with $k(v),n\geq0$. 
By Lemma~\ref{lem71}, if $v\overset{V}{\sim} w\in\Z^{m}$, then $k(v)=k(w)$. 
As a virtual analogue of Proposition~\ref{prop55}, 
we obtain the following. 

\begin{proposition}\label{prop79} 
For a nontrivial element $v\in\Z^{2m}$ $(m\geq2)$, 
the following are equivalent. 
\begin{enumerate}
\item
$v\overset{V}{\sim}
(\underbrace{x,\dots,x}_{\rm even}, 
\underbrace{y,\dots,y}_{\rm even})$ 
for some $x<y$. 
\item 
$\Delta(v)$ is divisible by $2^{k(v)+1}$. 
\end{enumerate} 
\end{proposition}

\begin{proof}
\underline{(i)$\Rightarrow$(ii).} 
By Lemmas~\ref{lem71} and~\ref{lem73}, we have 
\[\Delta(v)\equiv
\Delta(\underbrace{x,\dots,x}_{\rm even}, 
\underbrace{y,\dots,y}_{\rm even})
=0\pmod{2^{k(v)+1}}.\]

\underline{(ii)$\Rightarrow$(i).}
By Proposition~\ref{prop78}, 
we may assume that 
\[
v=(\underbrace{x,\dots,x}_{p}, 
\underbrace{y,\dots,y}_{2m-p})\in\Z^{2m}
\]
for some integers $x<y$ and $p\geq1$. 

Suppose, for contradiction, that $p$ is odd. 
A direct computation shows that 
\[
\Delta(v)=x-y 
\text{ and } 
d(v)=y-x.\] 
Since $y-x=2^{k(v)}(2n+1)$ for some $n\geq0$, 
we have $y-x\not\equiv0\pmod{2^{k(v)+1}}$, 
and hence 
\[
\Delta(v)\not\equiv 0\pmod{2^{k(v)+1}}.  
\]
This contradicts the assumption that $\Delta(v)\equiv0\pmod{2^{k(v)+1}}$. 
\end{proof}



\end{document}